\documentclass[10pt]{article}
 
\usepackage{hyperref}
\usepackage{amssymb, amsmath}
\input{liemacs10.sty}

\newcommand{\oS}{\oline{\bS}}
\newcommand{\oP}{\oline{\bP}}
\renewcommand{\mlabel}{\label}

\begin{document} 


\title{Positive energy representations and\\
 continuity of projective representations \\ 
for general topological groups} 
\author{Karl-Hermann Neeb
\begin{footnote}{
Department  Mathematik, FAU Erlangen-N\"urnberg, Cauerstrasse 11, 
91058 Erlangen, Germany, Email: karl-hermann.neeb@math.uni-erlangen.de} 
  \end{footnote}
\begin{footnote}{Supported by DFG-grant NE 413/7-2, Schwerpunktprogramm ``Darstellungstheorie''.}
\end{footnote}}

\maketitle


\begin{abstract} 
Let $G$ and $T$ be topological groups, 
$\alpha \: T \to \Aut(G)$ a homomorphism defining a continuous action 
of $T$ on $G$ and $G^\sharp := G \rtimes_\alpha T$ 
the corresponding semidirect product group. 
In this paper we address several issues concerning 
irreducible continuous unitary representations 
$(\pi^\sharp, \cH)$ of $G^\sharp$ whose restriction to 
$G$ remains irreducible. First we prove that, 
for $T = \R$, this is the case for 
any irreducible positive energy representation of $G^\sharp$, i.e., 
for which the one-parameter group $U_t := \pi^\sharp(\1,t)$ has non-negative spectrum.
The passage from irreducible unitary representations of $G$ to 
representations of $G^\sharp$ requires that 
certain projective unitary representations are continuous. 
To facilitate this verification, we derive various 
effective criteria for the continuity of projective unitary representations.
Based on results on Borchers for $W^*$-dynamical systems, 
we also derive a characterization of the continuous 
positive definite functions on $G$ that extend 
to a $G^\sharp$. \\
{\em Keywords:} positive energy representation, Borchers--Arveson Theorem, 
projective unitary representation, projective space \\
{\em MSC2010:} 22E45, 22E66. 22D10, 43A65  
\end{abstract}

\section*{Introduction} 

Let $G$ and $T$ be topological groups and 
$\alpha \: T \to \Aut(G)$ be a homomorphism defining a continuous 
action of $T$ on $G$. Then the 
semidirect product $G^\sharp := G \rtimes_\alpha T$ is also a 
topological group. In this note we take a closer look at irreducible 
continuous unitary representations of $G^\sharp$ whose restriction to 
$G$ remains irreducible. This is motivated by the concrete class of examples 
where $G$ is a Banach--Lie group, such as the group of $H^1$-maps 
$\bS^1 \to K$, where $K$ is a Lie group, and 
$\alpha \: \R \to \Aut(G)$ correspond to rotations of the circle.
Then $G^\sharp$ is not a Lie group, but it still is a topological group (cf.\ \cite{Ne12}).

For $T = \R$, a particularly interesting class of representations 
$(\pi^\sharp, \cH)$ of $G^\sharp$ 
are those for which the infinitesimal generator of the 
one-parameter group $U_t := \pi^\sharp(\1,t)$ has non-negative spectrum; 
the so-called {\it positive energy representations}. 
Our first main result (Theorem~\ref{thm:1.1}) asserts that every irreducible 
positive energy representation of $G^\sharp$ 
remains irreducible when restricted to $G$. This is quite remarkable 
because for non-trivial actions, there are many situations where 
irreducible representations of $G^\sharp$ are not irreducible 
when restricted to $G$. The simplest example are the 
irreducible representations of the orientation preserving 
affine group of $\R$, which is of the form 
$G^\sharp = \R \rtimes_\alpha \R$ for $\alpha_t(x) = e^t x$. 
We derive Theorem~\ref{thm:1.1}  as a consequence of the Borchers--Arveson 
Theorem on covariant representations of von Neumann algebras 
(cf.\ \cite[Thm.~3.2.46]{BR02}).

In the light of Theorem~\ref{thm:1.1} it becomes a 
natural question which irreducible unitary representations 
of $G$ extend to a positive energy representation of 
$G^\sharp = G \rtimes_\alpha \R$. 
This appears to be a difficult problem, and already for the 
rather concrete class of loop groups the answers becomes rather 
involved because it turns out that, for 
$G = C^\infty(\T,K)$ and $(\alpha_{z}f)(w) := f(zw)$ (the rotation action),  
smooth positive every representations of $G^\sharp$ are trivial in the sense that 
they factor through representations of products 
$A \times \R$, where $A$ is an abelian quotient of $K$. 
This is shown in Section~\ref{sec:5}, where we also derive a similar 
result for groups of the form 
$G = C^\infty(\R,K)$ and $(\alpha_tf)(s) = f(s+t)$. 
It is well-known from the classical theory of loop groups 
that the problem of the triviality of positive energy 
representations of $C^\infty(\T,K) \rtimes \T$ can be resolved 
by passing from the loop group $C^\infty(\T,K)$ to a suitable 
central extension $\tilde G$ for which 
$\tilde G \rtimes\T$ has plenty of positive energy representations 
 (\cite{PS86}, \cite{Ne12}). One obtains a similar picture 
for the case where $G$ is abelian, where one has to pass to Heisenberg groups 
to obtain positive energy representations (see \cite{NZ11} and \cite{Ze12} for details). 

Given an irreducible unitary representation 
$(\pi, \cH)$ of $G$, then a necessary condition for the 
existence of an extension to $G^\sharp = G \rtimes_\alpha T$ is that, for every 
$t \in T$, the representation $\pi \circ \alpha_t$ is equivalent 
to $\pi$. Then there exists, up to a phase factor in $\T$ 
uniquely determined, unitary operators $(U_t)_{t \in T}$ with 
\begin{equation}
  \label{eq:covar1}
 U_t \pi(g) U_t^* = \pi(\alpha_t(g)) \quad \mbox{ for } \quad t \in T, g \in G,
\end{equation}
and the question is whether these operators can be chosen in such a way 
that $U \: T \to \U(\cH)$ is a continuous group homomorphism. 
For any group homomorphism $\chi \: T \to \T$, the operators 
$\tilde U_t := U_t \chi(t)$ also satisfy the relation 
\eqref{eq:covar1} and since groups such as $T = \R$ have many discontinuous 
homomorphisms $\chi \: T \to \T$, one cannot hope for an arbitrarily 
chosen family $(U_t)$ to be continuous. However, what 
is uniquely determined by $\pi$ is the corresponding 
projective representation $\oline U \: T \to \PU(\cH)$, 
so that the main issue is whether this 
homomorphism is continuous or not. 

To facilitate the verification of the continuity of a 
projective unitary representation of a topological group, 
we provide in Section~\ref{sec:4} some useful criterion 
(Theorem~\ref{thm:3.5}). It asserts that, for a connected topological 
group $G$ and a projective unitary representation 
$\pi \: G \to \PU(\cH)$ continuity already follows from the existence 
of a cyclic ray $[v] = \C v$  in the projective space $\bP(\cH)$ of 
$\cH$ whose orbit map $G \to \bP(\cH), g \mapsto \pi(g)[v]$ is continuous. 
In particular, for irreducible representations, it suffices that 
one ray has a continuous orbit map. 
The main point of this criterion is that it is easier to check continuity 
for one vector than for all of them. 
This result has an analog for unitary representations $\pi \: G \to \U(\cH)$ 
whose proof is much 
easier, namely that a unitary representation generated by a vector 
with continuous orbit map is continuous. What makes the projective 
version more difficult is the lack of ``addition'' on projective space. 
We overcome this problem by using the midpoint 
operation for the natural Riemannian metric on $\bP(\cH)$ as a 
replacement. Our Theorem~\ref{thm:3.5} refines the criterion by Rieffel (\cite{Ri79}) 
that a projective unitary representation is continuous if and only if all 
its orbit maps in $\bP(\cH)$ are continuous.

Since the continuity of a projective unitary representation is 
characterized in terms of continuity of orbit maps in $\bP(\cH)$, 
it is also important to have effective tools to verify such continuity properties. 
In Section~\ref{sec:3} we provide such a tool 
by giving sufficient conditions for a subset $E\subeq \cH$ such that 
the topology on the projective space $\bP(\cH)$ is initial with 
respect to the functions $h_v([w]) := |\la v,w\ra|$, $w \in E$. 
In particular, it suffices that the corresponding functions $\ell_v(A) := \la Av,v\ra$ 
separate the points on the space of hermitian trace class operators. 
For any such set the continuity of an orbit map of 
a ray $[w] \in \bP(\cH)$ for a projective unitary 
representation can now be verified in terms of the continuity 
of the scalar-valued functions $G \to \R, g \mapsto h_v(\pi(g)[w])$. 

The homomorphism $\alpha \: T \to \Aut(G)$ 
is a group theoretic variant of the structure of a $C^*$-dynamical system, 
where $\alpha \: T \to \Aut(\cA)$ is a homomorphism into the group 
of automorphisms of a $C^*$-algebra $\cA$. 
For covariant representations $(\pi, U)$ of a $C^*$-algebra $\cA$ with respect to 
$\alpha$, Borchers gives in \cite{Bo83} a characterization of those states 
in $\cA^*$ (the topological dual of $\cA$) 
occurring for covariant representations (for which $U$ is continuous) 
as those transforming continuously under the natural $T$-action on the Banach space 
$\cA^*$. For covariant representations 
satisfying the spectral condition, i.e., the one-parameter group 
$(U_t)_{t \in \R}$ has non-negative spectrum, extra conditions 
on the states have to be imposed (cf.\ \cite[Sect.~II.5]{Bo96}). 
In Section~\ref{sec:6} we explain how Borchers result can be applied 
to the problem to extend continuous representations 
of $G$ to continuous representations 
of $G^\sharp = G \rtimes_\alpha T$. 
For the special case of pure states, this follows easily 
from our continuity criteria for projective representations.

For representations on separable Hilbert spaces and 
separable topological groups which are completely metrizable, so-called polish 
groups, one can also derive the continuity of projective unitary representations 
from rather weak measurability requirements. This has been shown 
by Cattaneo in \cite{Ca76} based on the results that,
for a separable Hilbert space $\cH$, its projective space 
$\bP(\cH)$ is a polish space (\cite[Prop.~4]{Ca76}) and 
$\PU(\cH)$ is a polish group (\cite[Prop.~5]{Ca76}).

{\bf Notation:} $\cH$ denotes a complex Hilbert space, 
$B(\cH)$ the $C^*$-algebra of bounded operators on $\cH$ 
and $\U(\cH)$ its unitary group. 
The {\it strong topology} on $\U(\cH)$ is the coarsest topology 
for which all functions 
$f_v(U) := \la Uv,v\ra$, $v \in \cH$, are continuous. 
It defines on $\U(\cH)$ a group topology and we write 
$\U(\cH)_s$ for this topological group. The center 
$\T \1 = Z(\U(\cH))$ is a closed subgroup, so that 
we also obtain the structure of a (Hausdorff topological group) 
on the projective unitary group $\PU(\cH) = \U(\cH)_s/\T \1$. 

The space of trace class operators on $\cH$ is denoted $B_1(\cH)$,
and we write $\Herm_1(\cH)$ for the subspace of hermitian trace class operators. 
For $v,w \in \cH$ we write $P_{v,w}(x) := \la x,w \ra v$ and 
$P_v := P_{v,v}$ for the corresponding rank-one operators on $\cH$. 

\tableofcontents

\section{Irreducible positive energy representations} 
\mlabel{sec:2}

Let $G$ and $T$ be topological groups and 
$\alpha \: T \to \Aut(G), t \mapsto \alpha_t$, 
be a homomorphism defining a continuous 
$T$-action on $G$. We write 
$ G^\sharp := G \rtimes_\alpha T$ 
for the corresponding semidirect product. 

\begin{defn} We call a pair $(\pi, U)$ of a continuous unitary representation 
$(\pi,\cH)$ of $G$ and a continuous unitary representation 
$(U,\cH)$ of $T$ {\it covariant} if 
\[ U_t \pi(g) U_t^* = \pi(\alpha_t g)\quad \mbox{ for } \quad g \in G, t \in T.\] 
This is equivalent to 
$\pi^\sharp(g,t) :=\pi(g) U_t$ to define a 
continuous unitary representation of the semidirect product 
$G^\sharp$. For $T = \R$, a covariant representation 
$(\pi, U)$ is said to be a 
{\it positive energy representation}\begin{footnote}
{In the context of covariant representations of operator algebras 
this condition is simply called the {\it spectral condition} or 
the {\it spectrum condition}.
}\end{footnote}
if  
the the infinitesimal generator $A = -i\frac{d}{dt}|_{t= 0}  U_t$ of the unitary one-parameter group 
$(U_t)_{t\in \R}$ has non-negative spectrum. We then also say that the unitary 
one-parameter group $(U_t)_{t \in \R}$ has {\it non-negative 
spectrum}. 
\end{defn}

The following theorem is a central result on covariant representations 
of operator algebras. 

\begin{thm}{\rm(Borchers--Arveson Theorem; \cite[Thm.~3.2.46]{BR02})} \mlabel{thm:1.1}
Let $(\alpha_t)_{t \in \R}$ be a $\sigma$-weakly continuous 
one-parameter group of automorphisms of a von Neumann algebra 
$\cM \subeq B(\cH)$, i.e., for each $\beta \in \cM_*$ (the predual of $\cM$)
 and $M \in \cM$, the 
function $t \mapsto \beta(\alpha_tM)$ is continuous. Then 
the following are equivalent:  
\begin{description}
\item[\rm(i)] There exists a strongly continuous unitary 
one-parameter group $(U_t)_{t \in\R}$ in $\U(\cH)$ with 
non-negative spectrum such that 
\[ \alpha_t(M) = U_t M U_t^*\quad \mbox{ for } \quad t \in \R, M \in \cM.\] 
\item[\rm(ii)] There exists a strongly continuous 
unitary one-parameter group $(U_t)_{t \in\R}$ in $\U(\cM)$ with 
non-negative spectrum such that 
\[ \alpha_t(M) = U_t M U_t^*\quad \mbox{ for } \quad t \in \R, M \in \cM.\] 
\end{description}
\end{thm}

\begin{rem} \mlabel{rem:1.2} If $(U_t)_{t \in \R}$ is a 
strongly continuous unitary one-parameter group, 
$A \in B_1(\cH) \cong B(\cH)_*$ is a trace class operator and $B \in B(\cH)$, then the 
function 
\[ t \mapsto \tr(AU_tBU_t^*) = \tr(U_t^*AU_tB) \] 
is continuous because the action of $\R$ on 
$B_1(\cH)$ defined by $(t,A) \mapsto U_t A U_t^*$ 
is strongly continuous. Therefore 
$\alpha_t(B) := U_tBU_t^*$ defines a $\sigma$-weakly continuous 
one-parameter group of automorphisms of $B(\cH)$, and hence on 
every von Neumann algebra $\cM \subeq B(\cH)$ which is invariant under conjugation 
with the operators $U_t$. Therefore  the $\sigma$-weak continuity of 
$\alpha$ is necessary for the conclusion of the Borchers--Arveson Theorem to hold. 
\end{rem}

\begin{cor}
  \mlabel{cor:1.3} 
If $(U_t)_{t \in \R}$ is a 
strongly continuous one-parameter group in $\U(\cH)$ 
with non-negative 
spectrum and $\cM \subeq B(\cH)$ be a von-Neumann algebra 
invariant under the automorphisms 
$\alpha_t(M) := U_t M U_t^*$ of $B(\cH)$, then there 
exists a strongly continuous unitary one-parameter group 
$(V_t)_{t \in \R}$ of $\cM$ with non-negative spectrum and a 
strongly continuous unitary one-parameter group 
$(W_t)_{t \in \R}$ in the commutant $\cM'$ such that 
\[ U_t = V_t W_t \quad \mbox{ for } \quad t \in \R.\] 
\end{cor}

\begin{prf} From Theorem~\ref{thm:1.1} we derive the existence of 
$(V_t)$ in $\U(\cM)$ satisfying
\[ U_t M U_t^* = \alpha_t(M) = V_t M V_t^* \quad \mbox{ for } \quad 
M \in \cM, t \in \R.\] 
Applying this relation to $V_s \in \cM$, we see that $U_t V_s = V_s U_t$ holds 
for $s,t \in \R$. Therefore $W_t := V_t^* U_t \in \U(\cM')$ defines a one-parameter group. 
\end{prf}

\begin{thm} \mlabel{thm:1.5} {\rm(Irreducibility Theorem)} 
If $(\pi^\sharp,\cH)$ is an irreducible positive energy 
representation of $G^\sharp = G \rtimes_\alpha \R$, 
then its restriction $\pi := \pi^\sharp\res_G$ is also irreducible. 
\end{thm}

\begin{prf} Let $\cM  := \pi(G)''$ be the von Neumann algebra 
generated by $\pi(G)$. Since $G \trile G^\sharp$ is normal, 
$\cM$ is invariant under conjugation with the unitary one-parameter group 
$U_t := \pi^\sharp(\1,t)$, so that Corollary~\ref{cor:1.3} 
provides a factorization 
$U_t = V_t W_t$ with strongly continuous one-parameter groups 
$(V_t)$ in $\cM$ and $(W_t)$ in $\cM' = \pi(G)'$. 
It follows in particular that 
$U_t \in \cM \cdot W_t$, so that 
\[ \pi^\sharp(G^\sharp)' = \pi(G)' \cap U_\R' = \cM' \cap U_\R'\supeq W_\R.\] 

If $\pi^\sharp$ is irreducible, 
$W_t \in \pi^\sharp(G^\sharp)' = \C \1$ for $t \in \R$. This implies 
that $U_t = V_t M_t \in \cM$, which finally leads to 
$\cM = \pi^\sharp(G^\sharp)'' = B(\cH)$, i.e., 
$\pi(G)' = \cM' = \C \1$, which means that $\pi$ is irreducible. 
\end{prf}

\begin{rem} \mlabel{rem:1.6} 
If $(\pi,\cH)$ is an irreducible representation 
of $G$ which has some extension $\pi^\sharp$ to $G^\sharp 
= G \rtimes_\alpha \R$, then 
Schur's Lemma implies that this extension is unique up to a unitary character of $\R$. 
For each $\mu\in \R$, we obtain a modified unitary representation 
\[ \pi^\sharp_\mu(g,t) = e^{it\mu}\pi^\sharp(g,t).\] 
The set of all those $\mu \in \R$ 
for which $\pi^\sharp_\mu$ is a positive energy representation is an interval 
which is either empty or of the form $[\mu_0, \infty[$ for some $\mu_0 \in \R$. 
In the latter case $\mu_0 = - \inf \Spec(-i\dd\pi^\sharp(0,1))$ leads to a representation 
for which $\inf\Spec(-i\dd\pi^\sharp_{\mu_0}(0,1)) = 0$. We call this the 
{\it minimal positive energy extension} of $\pi$. 

Since the group $G^\sharp$ has richer structure than $G$ itself, the 
crucial advantage of Theorem~\ref{thm:1.5} is that it permits us to study 
certain irreducible representations of $G$ as representations of $G^\sharp$ 
(cf.\ \cite{Ne12}, \cite{PS86}). 
\end{rem}

\section{Some Lie group examples} 
\mlabel{sec:5} 

In this section $G$ denotes a Lie group modeled on a locally convex space 
(see \cite{Ne06} for more details). 
If $G$ is a Lie group and the action $\alpha$ of $\R$ on $G$ is smooth, the group 
$G^\sharp$ also is a Lie group, so that it makes sense to consider unitary representations $(\pi^\sharp, \cH)$ of 
$G^\sharp$ which are smooth in the sense that the subspace $\cH^\infty$ of smooth vectors 
is dense in $\cH$. The following lemma, 
applied with $y = (0,1) \in \g^\sharp = \g \rtimes \R$, 
the Lie algebra of $G^\sharp$, 
shows that, in some situations 
the positive energy condition leads to serious restrictions on the 
corresponding smooth representation $(\pi, \cH)$. 

\begin{lem} \mlabel{lem:solvlem} {\rm(\cite[Lemma~5.12]{Ne12})} Let $(\pi, \cH)$ be a smooth representation of a 
Lie group $G$ with exponential function $\exp \: \g \to G$ 
and $x,y \in \g$ with $[x,[x,y]]~=~0$. 
If $-i\dd\pi(y)$ is bounded from below, then $\dd\pi([x,y]) = 0$. 
\end{lem}

With the preceding lemma one can show in particular that: 

\begin{thm} Let $G$ be a connected Lie group 
with Lie algebra $\g = C^\infty(\T,\fk)$, $\T = \R/\Z$ and $K$ a connected Lie group 
with Lie algebra $\fk$. Writing $[t] := t + \Z$ for elements of $\T$ we obtain an action 
of $\T$ on $G$ by 
\[ (\alpha_{[t]}f)([s]) = f([s-t]).\] 
Then all smooth positive energy representation of 
$G^\sharp = G \rtimes_\alpha \R$ factor through 
representations of an abelian quotient group of~$K$ 
and annihilate the commutator group of $G$. 
Conversely, all such representations are positive energy 
representations for trivial reasons.  
\end{thm}

\begin{prf} (a) First we consider the case $G = C^\infty(\T, \R)$. 
Applying Lemma~\ref{lem:solvlem} with $y = (0,1)$ and $x \in \g$, it follows that 
$x' = [y,x] \in \ker(\dd\pi)$, so that 
\[ \Big\{ \xi \in \g = C^\infty(\T,\R) \: \int_\T  \xi(t)\, dt = 0\Big\} \subeq \ker(\dd\pi).\] 

(b) For a general Lie group $K$, we obtain for each $x \in \fk$ a smooth homomorphism 
\[ \Gamma_x \: C^\infty(\T,\R) \to C^\infty(\T,K), \quad 
\Gamma_x(f)(t) := \exp(f(t)x).\] 
As $C^\infty(\T,\R) \otimes \fk$ is dense in $\g$, we obtain from (a) the relation 
\[ \Big\{ \xi \in \g = C^\infty(\T,\fk) \: \int_\T  \xi(t)\, dt = 0\Big\} \subeq \ker(\dd\pi).\] 
For $\T = \R/2\pi \Z$, 
$\xi_1(t) = \cos t \cdot x$, 
$\xi_2(t) = \sin t \cdot x$, 
$\xi_3(t) = \cos t \cdot y$ and 
$\xi_4(t) = \sin t \cdot y$, the relation 
\[ [\xi_1, \xi_3] + [\xi_2, \xi_4]  = [x,y]\]  
now shows that $[\fk,\fk] \subeq \ker(\dd\pi)$, which leads to 
\[ \Big\{ \xi \in \g 
= C^\infty(\T,\fk) \: \int_\T  \xi(t)\, dt \in \oline{[\fk,\fk]}\Big\} \subeq \ker(\dd\pi).\] 
Therefore the derived representation of $\g$ actually factors through a representation of the abelian 
quotient algebra $\fk/\oline{[\fk,\fk]} \cong \g/\oline{[\g,\g]}$. 
\end{prf}

For the following theorem we recall that a Lie group 
$K$ is called {\it regular} if, for each smooth map
$\xi \: [0,1] \to \fk$, the Lie algebra of $K$, there exists a smooth curve 
$\gamma_\xi \: [0,1] \to K$ with 
$\gamma_\xi(0) = \1$ and $\gamma_\xi'(t) = \gamma_\xi(t) \cdot \xi'(t)$ 
(here $\cdot$ refers to the canonical action $K \times TK \to TK$)  
and the map $\evol \: C^\infty([0,1],\fk) \to G, \xi \mapsto \gamma_\xi(1)$ is smooth. 
For any regular Lie group $K$, the parametrization of smooth curves
$\gamma \: \R \to K$ by their logarithmic derivatives $\xi \: \R \to \fk$ 
leads to a natural Lie group structure on the group 
$G := C^\infty(\R,K)$, endowed with the pointwise multiplication 
(see \cite{NW08} for details). 

\begin{thm} Let $K$ be a connected regular Lie group 
with Lie algebra $\fk$ and 
$G := C^\infty(\R,K)$. Then $(\alpha_t f)(s) = f(s-t)$ defines a smooth action of $\R$ on $G$. 
Then all smooth positive energy representation of 
$G^\sharp = G \rtimes_\alpha \R$ are trivial on~$G$. 
\end{thm}

\begin{prf} The infinitesimal generator of the derived action 
of $\R$ on $\g$ is 
$Df := \frac{d}{dt}|_{t = 0} \L(\alpha_t)f = - f'$. 
For $y = (0,1)$ and $f(\R) \subeq \R x$ for some $x \in \fk$, 
we derive from Lemma~\ref{lem:solvlem} that $\dd\pi(f) = 0$ for every positive energy 
representation of $G^\sharp$. Therefore 
$C^\infty(\R) \otimes \fk \subeq \ker \dd\pi$, and since this subspace of 
$\g$ is dense, we obtain $\dd \pi(\g)~=~0$. 
\end{prf}

\section{Generating the topology on projective space} 
\mlabel{sec:3} 

In this section we discuss  the topology 
of the projective space $\bP(\cH)$. Our main goal 
is a criterion for a subset $E\subeq \cH$ such that 
the topology on $\bP(\cH)$ is initial with respect to the 
functions $h_v([w]) := |\la v,w\ra|$, $v \in E$. 
It suffices that the corresponding functions $\ell_v(A) := \la Av,v\ra$ 
separate the points on the space of hermitian trace class operators. 
Here we use the topological embedding 
$\bP(\cH) \into \Herm_1(\cH), [w] \mapsto P_w$ of $\bP(\cH)$ 
in terms of rank-one projections into the space $C$ of positive 
trace class operators $A$ with $\tr(A) \leq 1$. 
It is crucial to our argument that, if $\cH$ is infinite-dimensional,
then $C$ can be considered as a compactification of $\bP(\cH)$, so that 
compactness arguments can be used.

\begin{defn} Let $\cH$ be a complex Hilbert space, endow its 
unit sphere $\bS(\cH) := \{v \in \cH \: \|v\| =1\}$ with the subspace topology 
inherited from $\cH$ and the projective space 
$\bP(\cH) \cong \bS(\cH)/\T$ with the quotient topology. We denote its elements, the 
one-dimensional subspaces of $\cH$, by $[v] = \C v$. 
\end{defn}

Since the sphere inherits a natural metric from $\cH$, it is instructive to first 
take a closer look on the metric aspects of the topology on $\bP(\cH)$: 

\begin{lem} \mlabel{lem:2.1} 
{\rm (a)} The metric  $d(x,y) = \|x-y\|$ on the sphere $\bS(\cH)$ 
induces on $\bP(\cH)$ the metric 
\[ d([x],[y]) := d(\T x,\T y) = \sqrt{2(1- |\la x,y\ra|)} \in [0,\sqrt 2], \quad 
x,y \in \bS(\cH).\] 

{\rm(b)} The map $\iota  \: \bP(\cH) \into \Herm_1(\cH), [v] \mapsto P_v$, 
is a topological embedding. 
\end{lem}

\begin{prf} (a) follows from
$d(\T x,\T y) = d(x,\T y) = \inf_{|t|=1} d(x,ty)$ and 
\[ \inf_{|t|=1} d(x,ty)^2 
= \inf_{|t|=1} \|x-ty\|^2 
= \inf_{|t|=1} 2(1-\Re \la x,ty\ra)
= 2(1- |\la x,y\ra|).\] 

(b) For $v,w \in \bS(\cH)$ with $[v]\not=[w]$ 
the operator $A := P_v - P_w$ is hermitian of rank $2$  
with $\tr(A) = \|v\|^2 - \|w\|^2 = 0$. 
Write $w = \lambda v + \mu v'$ with $v' \in \bS(\cH)$ orthogonal 
to $v$. Then 
\[ P_w = |\lambda|^2 P_v + |\mu|^2 P_{v'} 
+ \lambda\oline{\mu} P_{v,v'} + \oline\lambda\mu P_{v',v}.\] 
On the subspace $\cH_0 = \C v + \C w$ we therefore have 
\[ \det(A\res_{\cH_0}) 
= (1- |\lambda|^2)(-|\mu|^2) - |\lambda|^2|\mu|^2 
= -|\mu|^2.\] 
As $\tr(A\res_{\cH_0}) = 0$, the two eigenvalues of $A\res_{\cH_0}$ are 
$\pm |\mu|$, which leads to 
\[ \|P_v - P_w\|_1 = 2|\mu| = 2\sqrt{1 - |\lambda|^2} 
= 2\sqrt{1 -|\la v,w\ra|^2}.\] 
This implies the assertion. 
\end{prf}

Now we turn from the metric point of view to families of functions generating 
the topology. 

\begin{defn} \mlabel{def:fvhv} For $v\in \cH$ we define the functions 
\[ e_v \: \cH \to \C, \quad e_v(x) := \la v, x \ra\]
and 
\[  h_v \: \cH/\T  \to \C, \quad h_v([x]) := |\la v, x \ra| \quad \mbox{ for } \quad 
[x] = \T x.\] 
\end{defn}

\begin{rem} \mlabel{rem:3.1} 
(a) The topology on the sphere $\bS(\cH)$ 
is the initial topology defined by the functions 
$e_v$, $\|v\| =1$, because it coincides with the weak 
topology. In fact, for $x, y \in \bS(\cH)$, we have
\[ \|x-y\|^2 = 2(1 - \Re \la x,y\ra).\] 

Let  $E \subeq \cH$ is a total subset, i.e., $\Spann E$ is dense in $\cH$. We 
consider on $\bS(\cH)$ the initial topology $\tau_E$ defined by the functions 
$\{e_v \: v \in E\}$. Then 
\[ \{ w \in \cH \: e_w \in C(\bS(\cH),\tau_E)\} \] 
is a closed subspace of $\cH$ 
because $\|e_w\|_\infty \leq \|w\|$. Since this subspace contains $E$, it coincides
with~$\cH$. Therefore $\tau_E$ coincides with the metric topology on~$\bS(\cH)$. 

(b) From the formula for the quotient metric in Lemma~\ref{lem:2.1}(a), it follows immediately 
that the topology on $\bP(\cH)$ is the initial topology defined by the functions 
$h_v$, $\|v\| = 1$. 
\end{rem} 

\begin{rem} \mlabel{rem:3.5n}
It is an interesting question which conditions we have
to require from a subset $E \subeq \cH$ to ensure that 
$h_E := \{ h_v \: v \in E\}$ defines the topology on 
$\bP(\cH)$. It is certainly necessary that $h_E$ separates the points of 
$\bP(\cH)$. For any proper orthogonal decomposition 
$\cH = \cH_1 \oplus \cH_2$, this rules out subsets $E \not\subeq \cH_1 \cup \cH_2$ 
because then $h_E$ cannot separate the elements 
$[v_1 + \zeta v_2] \in \bP(\cH)$, $\zeta \in \T$, where 
$0 \not= v_j \in \cH_j$, $j =1,2$. 
This implies in particular that (for $\dim \cH > 1$)  $E$ needs to be 
total and that it cannot be decomposed into two proper 
mutually orthogonal subsets. However, this 
condition is not sufficient for $h_E$ to separate the points of 
$\bP(\cH)$. For $\cH = \C^2$ and 
a non-orthogonal basis $v_1, v_2 \in \C^2$, the functions 
$h_{v_1}$ and $h_{v_2}$ do not separate the points of $\bP(\cH)$. 
The level sets of both functions are families of circles on 
the Riemann sphere $\bP(\cH) \cong \bS^2$ and two such circles can intersect in 
two points. Geometrically this means that a ray 
$\C v \subeq \C^2$ is not determined by the two numbers 
$|\la v, v_1\ra|$ and $|\la v, v_2\ra|$. 
From the same reasoning it follows that, if 
$v_1, v_2, v_3 \in \C^2$ are such that the corresponding rays 
$[v_1], [v_2], [v_3] \in \bP(\C^2) \cong \bS^2$ do not lie on any great circle, 
then the functions $h_{v_j}$, $j=1,2,3$, separate the points of $\bP(\C^2)$.   
\end{rem}

\begin{prob} Suppose that $h_E$ separates the points of $\bP(\cH)$. Is the topology on 
$\bP(\cH)$ the initial topology with respect to the set $h_E$? 
\end{prob}

Below we prove a slightly weaker statement, which requires $h_E$ to separate the functions 
on a slightly larger set that we introduce below.  

\begin{defn} Let 
\[ \oS(\cH) := \{ v \in \cH \: \|v\| \leq 1\} \] 
denote the closed unit ball in $\cH$, endowed with the weak topology, with respect to which 
it is a compact space. Since the scalar multiplication action 
$\T \times \oS(\cH) \to \oS(\cH)$ is continuous and $\T$ is compact, we obtain on the 
quotient space 
\[ \oP(\cH) := \oS(\cH)/\T \] 
a compact Hausdorff topology and a topological embedding
$\bP(\cH) \into \oP(\cH).$ 
\end{defn}

\begin{rem} \mlabel{rem:2.6} (a) If $\cH$ is infinite-dimensional, then $\bS(\cH)$ is dense in 
$\oS(\cH)$, so that we may consider the ball as a compactification of the unit sphere. 
Likewise $\oP(\cH)$ is a compactification of $\bP(\cH)$. 
If $\cH$ is finite-dimensional, the unit sphere and the projective space are compact. 

(b) If $X$ is a compact (Hausdorff) space and 
$\cE \subeq C(X,\C)$ is a point separating set of continuous functions, then the topology on $X$  
coincides with the initial topology defined by $\cE$ because the map 
$X \to \C^\cE, x \mapsto (f(x))_{f \in \cE}$ is a topological embedding. 
\end{rem}

\begin{prop} \mlabel{prop:2.7} If $E \subeq \cH$ is a subset for which the functions 
$h_v([x]) := |\la v,x\ra|$, $v \in E$, separate the points of $\oP(\cH)$, then 
the topology on $\bP(\cH)$ is the initial topology with respect to the family 
$(h_v)_{v \in E}$. This is in particular the case 
for $E = \cH$. 
\end{prop}

\begin{prf} (a) Since the functions $x \mapsto |\la x,v\ra|$ on $\oS(\cH)$ are continuous, 
the function $h_v$ defines a continuous functions on $\oP(\cH)$. 
In view of Remark~\ref{rem:2.6}(b), the topology on the compact space 
$\oP(\cH)$ is initial with respect to the functions $(h_v)_{v \in E}$. 
Now the first assertion follows from the fact that $\bP(\cH)$ is a topological 
subspace of $\oP(\cH)$. 

(b) To verify the second assertion, we show that, 
for $E = \bS(\cH)$, the functions $(h_v)_{v \in E}$ separate the points of $\oP(\cH)$. 
In fact, they obviously separate $0$ from the non-zero elements in $\oS(\cH)$. 
Moreover, 
$\|x\| = \sup_{\|v\| = 1} h_v(x)$ 
implies that they determine the norm of an element. 
If $h_v(x) = h_v(y)$ holds for two non-zero elements 
$x,y \in \oline\bS(\cH)$ and all $v \in \bS(\cH)$, we obtain 
for $v = x/\|x\|$ the relation 
\[ \|x\| =  h_v(x) = h_v(y) = \frac{|\la x, y\ra |}{\|x\|} 
\quad \mbox{ and likewise } \quad  \|y\| = \frac{|\la x, y\ra| }{\|y\|}, \] 
so that 
\[ |\la x,y \ra| = \|x\| \cdot \|y\|.\] 
We conclude that $y \in \C x$, and since the preceding discussion also implies that 
$\|y\| = \|x\|$, it follows that $y \in \T x$, i.e., $[x] = [y]$ in $\oline\bP(\cH)$.   
\end{prf}

\begin{ex} \mlabel{ex:1} Not every family $(h_v)_{v \in E}$ which separate the points of 
$\bP(\cH)$ also separate the points of $\oP(\cH)$. A simple example arises for 
$\dim \cH = 1$, where $E = \eset$ suffices to separate the points of the one 
point set $\bP(\cH)$, but this is not enough for the interval 
$\oP(\cH)$.

For $\cH = \C^2$, we have seen in 
Remark~\ref{rem:3.5n} that if $v_1, v_2, v_3$ are three unit vectors for which the 
corresponding rays $[v_1], [v_2], [v_3] \in \bP(\C^2) \cong \bS^2$ do not lie on a great circle, 
the functions $h_{v_j}$, $j =1,2,3,$ separate the points of $\bP(\C^2)$. 
Since the space $\Herm_2(\C)$ of hermitian $2 \times 2$-matrices is 
$4$-dimensional, there exists a non-zero matrix $A \in \Herm_2(\C)$ with 
\[ \la Av_j, v_j \ra = 0 \quad \mbox{ for } \quad j=1,2,3.\] 
If $A$ is positive or negative semidefinite, then these relations imply 
$Av_j = 0$ for $j =1,2,3$, and since the $v_j$ are linearly independent, 
this contradicts $A \not=0$. Therefore 
$A$ has eigenvalues ${\lambda_1 < 0 < \lambda_2}$. We assume w.l.o.g.\ that 
$|\lambda_j| \leq 1$ for $j =1,2$ and write 
$u_1, u_2 \in \C^2$ for corresponding unit eigenvectors of $A$ with 
$Au_j = \lambda_j u_j$. From $A = \lambda_1 P_{u_1} + \lambda_2 P_{u_2}$ we 
then obtain for $w_j := \sqrt{|\lambda_j|} u_j$ the relation 
$A = P_{w_2} - P_{w_1}$, and thus 
\[ 0 = \la A v_j, v_j\ra = |\la w_2, v_j \ra|^2 - |\la w_1, v_j \ra|^2 \] 
implies that the functions $h_{v_j}$ do not separate the two elements 
$[w_1]$ and $[w_2]$ in $\oline\bP(\C^2)$. 
\end{ex}

We have already seen in Lemma~\ref{lem:2.1}(b) that we have a topological embedding 
$\eta \: \bP(\cH) \into \Herm_1(\cH), [v] \mapsto P_v$. 
The subset 
\[ C := \{ A = A^* \in \Herm_1(\cH) \: 0 \leq A, \tr A \leq 1\}\] 
of $\Herm_1(\cH)$ is convex, bounded and weak-$*$-closed if we consider 
$B_1(\cH)$ via the trace pairing as the dual space of the space $K(\cH)$ 
of compact operators on $\cH$. We conclude that 
$C$ is a weak-$*$-compact subset. Next we observe that $\eta$ extends to a map 
\[ \eta \: \oP(\cH) \into C, \quad [v] \mapsto P_v. \]
To see that this map is continuous, we first recall that the subset 
$\{ P_v \: v \in \cH\} \subeq K(\cH)$ spans a dense subspace, which implies that the 
topology on $C$ is the initial topology with respect to the functions 
\begin{equation}
  \label{eq:lv}
 \ell_v \: C \to \C, \quad A \mapsto \tr(A P_v) = \tr(P_{Av,v}) = \la Av,v\ra
\end{equation}
(Remark~\ref{rem:2.6}(b)). 
Therefore 
\[ \ell_v(P_w) = \la P_w v, v\ra = \la v,w\ra \la w, v\ra = |\la v,w\ra|^2\] 
shows that all function $\ell_v \circ \eta$ are continuous 
(cf.\ Definition~\ref{def:fvhv}), so that 
$\eta \: \oP(\cH) \to C$ is continuous. From $\|P_v\| = \|v\|^2$ it further follows that 
$\eta$ is injective, hence a topological embedding of $\oP(\cH)$ onto a weak-$*$-compact subset 
of $C$. 

This leads to the following criterion: 

\begin{prop} \mlabel{prop:2.9} 
If $E \subeq \cH$ is such that the functions $(\ell_v)_{v \in E}$ separate the 
points of $\Herm_1(\cH)$, then $(h_v)_{v \in E}$ separates the points of $\oP(\cH)$. 
In particular, 
the topology on $\bP(\cH)$ is the initial topology with respect to $(h_v)_{v \in E}$. 
\end{prop}

\begin{lem} \mlabel{lem:pointsep} 
For a subset $E \subeq \cH$ the functions $(\ell_v)_{v \in E}$ separate the points of $C$ if and only 
if $\{ P_v \: v \in E\}$ spans a dense subspace of $K(\cH)$.   
\end{lem}

\begin{prf} The set $\{ P_v \: v \in E \}$ is not total in $K(\cH)$, i.e., 
its closed span is a proper subspace, if and only if 
there exists a hermitian trace class operator $A \in \Herm_1(\cH) \cong K(\cH)'$ with 
$\la Av,v\ra = 0$ for every $v \in E$. 

Writing $A = A_+ - A_-$ with 
positive operators $A_\pm$, we find a $\lambda > 0$ such that 
$\tr(\lambda A_\pm) < 1$. Then 
$\lambda A_\pm \in C$ satisfy 
$\la A_+ v,v\ra = \la A_-v,v\ra$ for every $v \in E$. Therefore 
$(\ell_v)_{v \in E}$ does not separate the points of $C$. 

If, conversely, $A_\pm \in C$ are two different operators not separated by the functions 
$(\ell_v)_{v \in E}$,
then $A := A_+ - A_- \in \Herm_1(\cH)$ is non-zero with 
$\tr(A P_v) = 0$ for every $v \in E$. Therefore $\{ P_v \: v \in E \}$ is not total. 
\end{prf}

\begin{rem} If $\dim \cH = 2$, then $\dim \Herm_1(\cH)= \dim \Herm_2(\C) = 4$, so that 
any subset $E \subeq \cH$ for which $(\ell_v)_{v \in E}$ separate the points of 
$\Herm_1(\cH)$ has to contain at least $4$ elements 
(cf.\ Example~\ref{ex:1}). 
\end{rem}

The following criterion is sometimes useful to verify the condition 
in Lemma~\ref{lem:pointsep}. 

\begin{prop} If $M$ is a connected complex manifold and $F \: M \to \cH$ a holomorphic map with 
total range, then $(\ell_{F(m)})_{m \in M}$ separates the points of $\Herm_1(\cH)$. 
\end{prop}

\begin{prf} Let $\oline M$ denote the real manifold $M$, endowed with the opposite 
complex structure. Then, for each $A \in \Herm_1(\cH)$, the function 
$\alpha_A(m,n) := \la A F(m), F(n) \ra$  on $M \times \oline M$ is holomorphic. 
If $\ell_{F(m)}(A) = \la A F(m), F(m) \ra = 0$ for every $m \in M$, then 
the holomorphic function $\alpha_A$  vanishes on the totally real 
submanifold $\Delta_M = \{ (m,m) \: m \in M\}$ of $M \times \oline M$, 
which implies that $\alpha_A = 0$ 
(cf.\ \cite[Prop.~A.III.7]{Ne00}). We conclude that, for each $m \in M$, 
$A F(m) \in \im(F)^\bot = \{0\}$, and since $\im(F)$ is total, it follows that $A = 0$. 
\end{prf}

The preceding proposition applies in particular 
to all reproducing kernel Hilbert spaces of holomorphic functions 
or holomorphic sections (cf.\ \cite{PS86}, \cite{Ne11}). 

\section{Continuity of projective unitary representations} 
\mlabel{sec:4}

It is well-known that a unitary representation 
$\pi \: G \to \U(\cH)$ of a topological group 
is continuous if and only if, for each $v$ in a total subset 
$E \subeq \cH$, the function 
$\pi^{v,v}(g) := \la \pi(g)v,v\ra$ is continuous 
(cf.\ \cite[Lemma~VI.1.3]{Ne00}). In this 
section we discuss a similar continuity criterion for 
projective unitary representations $\pi \: G \to \PU(\cH)$.

\begin{prop} \mlabel{prop:topprogrp} {\rm(a)} 
The topology on $\PU(\cH)$ is the coarsest 
topology for which all functions 
$$ h_{v,w} \: \PU(\cH) \to \R, \quad [g] \mapsto |\la gv, w \ra|, \quad 
v,w \in \cH, $$ 
are continuous. 

{\rm(b)} The quotient map $q \: \U(\cH) \to \PU(\cH)$ 
has continuous local sections, i.e., each $[g] \in \PU(\cH)$ 
has an open neighborhood $U$ on which there exists a 
continuous section $\sigma \: U \to \U(\cH)$ of~$q$. 
\end{prop}

\begin{prf} (a) Let $q \: \U(\cH) \to \PU(\cH)$ denote the quotient map. 
Then all functions $f_{v,w} := h_{v,w} \circ q$ are continuous on 
$\U(\cH)$, which implies that the functions $h_{v,w}$ are continuous 
on $\PU(\cH)$. 

Let $\tau$ denote the coarsest topology on $\PU(\cH)$ for which 
all functions $h_{v,w}$ are continuous. We know already that this 
topology is coarser than the quotient topology. 
Next we observe that the relations 
$$ h_{v,w}([g][g']) = h_{g'v,w}([g]) = h_{v,g^{-1}w}([g']) $$
imply  that left and right multiplications are continuous 
in $\tau$. To see that $\tau$ coincides with the quotient topology, 
it therefore remains to see that 
$[g_i] \to \1$ in $\tau$ implies that $[g_i] \to \1$ in the quotient 
topology. 

For a net $([g_i])_{i \in I}$ in $\PU(\cH)$ we consider 
a lift $(g_i)_{i \in I}$ in $\U(\cH)$. Since 
the closed operator ball $\cB := \{ A \in B(\cH) \: \|A\| \leq 1 \}$ 
is compact in the weak operator topology, there exists a convergent subnet 
$g_{\alpha(j)} \to g_0 \in \cB$. For $v, w \in \cH$ we then have
$$ h_{v,w}(g_{\alpha(j)}) \to h_{v,w}(\1) = |\la v,w \ra| $$
and also 
$$  h_{v,w}(g_{\alpha(j)}) = |\la g_{\alpha(j)}v, w\ra| 
\to |\la g_0 v, w\ra|, $$ 
hence $|\la g_0 v, w\ra| = |\la v,w\ra|$. This implies in 
particular that, for each non-zero vector $v$, we have 
$$g_0 v \in (v^\bot)^\bot = \C v, $$
so that each vector is an eigenvector, and this implies that 
$g_0 = t \1$ for some $t \in \C$. If $v = w$ is a unit 
vector, we obtain $|t| = |\la g_0v,v\ra| = 1$. 
Therefore we have $g_{\alpha(j)} \to t \1$ in $\U(\cH)$, and this 
implies that $[g_{\alpha(j)}] \to [\1]$ in $\PU(\cH)$. 

If the net $(g_i)_{i \in I}$ does not converge to 
$\1$ in $\PU(\cH)$, then there exists an open $\1$-neighborhood 
$U$ for which the set $I_U := \{ i \in I \: g_i \not\in U\}$ 
is cofinal, which leads to a subnet $(g_i)_{i \in I_U}$ 
converging to $\1$ in $\tau$ and contained in the closed subset $U^c$. 
Applying the preceding argument to this subnet now leads to a 
contradiction since it cannot have any subnet converging to $\1$ 
because $U^c$ is closed. 

(b) Since we can move sections with left multiplication maps, 
it suffices to assume that $g = \1$. 
Pick $0 \not= v_0 \in \cH$. Then 
$$ \Omega := \{ g \in \U(\cH) \: \la gv_0, v_0 \ra \not= 0\} $$
is an open $\1$-neighborhood in $\U(\cH)_s$ with 
$\Omega \T = \Omega$. Therefore $\tilde\Omega := \{ [g] \: g \in \Omega\}$ 
is an open $\1$-neighborhood of $\PU(\cH)$. For each 
$g \in \Omega$ there exists a unique $t \in \T$ with 
$$ tg \in \Omega_+ :=  \{ g \in \U(\cH) \: \la gv_0, v_0 \ra >0\}. $$
We now define a map 
$$ \sigma \: \tilde\Omega \to \Omega, \quad 
[g] \mapsto g \quad \mbox{ for } \quad g \in \Omega_+. $$
To see that $\sigma$ is continuous, it suffices to observe 
that the map 
$$ \Omega \to \Omega_+, \quad g 
\mapsto \frac{|\la g v_0, v_0\ra|}{\la g v_0, v_0\ra} g $$
is continuous and constant on the cosets of $\T$. Hence 
it factors through a continuous map $\tilde\Omega \to \Omega_+$ 
which is $\sigma$. 
Therefore the quotient map 
$$ q \: \U(\cH) \to \PU(\cH), \quad g \mapsto [g] $$
has a continuous section in the $\1$-neighborhood $\tilde\Omega$ of $\PU(\cH)$.
\end{prf}

\begin{cor} \mlabel{cor:3.2} Let $G$ be a topological group and 
$\pi \:G \to \PU(\cH)$ be a group homomorphism. Then the following are equivalent: 
\begin{description}
\item[\rm(i)] $\pi$ is continuous. 
\item[\rm(ii)] For all $v,w \in \cH$, the function $G \to \R, \mapsto |\la \pi(g)v,w\ra|$ 
is  continuous. 
\item[\rm(iii)] For each $[v] \in \bP(\cH)$, the orbit map 
$G \to \bP(\cH), g \mapsto \pi(g)[v]$ is continuous. 
\end{description}
\end{cor}

\begin{prf} The equivalence of (i) and (ii) is an immediate consequence of 
Proposition~\ref{prop:topprogrp}. Clearly, (ii) is equivalent to the requirement, that for 
any $0 \not= v \in \cH$, the functions 
$G \to h_w(\pi(g)[v])$, $w \in \cH$, are continuous (cf.\ Definition~\ref{def:fvhv}). 
As the topology on $\bP(\cH)$ is initial with respect to the 
functions $(h_w)_{w \in \cH}$ (Proposition~\ref{prop:2.7}), the corollary follows. 
\end{prf}

The equivalence of (i) and (iii) in the preceding corollary can already 
be found in \cite[Lemma~8.1]{Ri79} (see also \cite[Prop.~6]{Ca76}). 
The characterization of the continuity of a projective 
representation in Corollary~\ref{cor:3.2} involves all elements 
$v \in \cH$. This makes it inconvenient to use in practice. 
Below we develop a criterion which makes it much easier to check continuity.

\begin{defn} Let $G$ be a topological group. For a homomorphism 
$\pi \:G \to \PU(\cH)$ we write 
$\bP(\cH)_c \subeq \bP(\cH)$ for the set of all elements 
$[v]$ for which the $G$-orbit map is continuous. 
\end{defn}

We now take a closer look at the structure of the set $\bP(\cH)_c$. We start 
with an elementary observation on isometric actions of topological groups. 

\begin{lem} \mlabel{lem:closed} If the topological group $G$ acts isometrically by 
\[ \sigma \: G \times X \to X, \quad (g,x) \mapsto g.x \] 
 on the metric space $(X,d)$, then 
the set $X_c$ of all points with continuous orbit maps is closed and 
the $G$-action on this set is continuous. 
\end{lem}

\begin{prf} If $x_n \to x$ in $X$, then the orbits maps 
$\sigma^{x_n} \: G \to X$ converge uniformly to the orbit map 
$\sigma^x$. This proves that $X_c$ is closed. 

The second assertion follows from 
\[ d(g_0.x_0, g_1.x_1)  = d(g_1^{-1}g_0.x_0, x_1) 
\leq d(g_1^{-1}g_0.x_0, x_0) + d(x_0, x_1).
\qedhere\]   
\end{prf}

\begin{rem} For a unitary representation $\pi \: G \to \U(\cH)$, one can also 
consider the subset $\cH_c \subeq \cH$ of those elements with a 
continuous orbit map. This is obviously a subspace which is closed by 
Lemma~\ref{lem:closed}. Now $\cH = \cH_c \oplus \cH_c^\bot$ 
provides a decomposition of the unitary representation as a direct 
sum of a continuous representation and one without non-zero continuous 
orbit maps. We shall see below that the situation is more complicated in 
the projective case.   
\end{rem}

\begin{defn}
On $\bP(\cH)$ we now consider the Riemannian metric given by 
\[ d_R([x],[y]) := \arccos |\la x,y\ra| \in [0,\pi/2] \quad 
\mbox{ for } \quad x,y \in \bS(\cH).\] 
We write $[x] \bot [y]$ and say that $[x]$ and $[y]$ are 
{\it orthogonal} if $x \bot y$, which is equivalent to 
$d_R([x],[y]) = \pi/2$. 
If $[x]$ and $[y]$ are not orthogonal, then we write 
$[x] \sharp [y]$ for the unique metric midpoint of $[x]$ and $[y]$. 
\end{defn}

\begin{lem} \mlabel{lem:mipo} The midpoint operation is a continuous 
map 
\[ \sharp \: \{([x], [y]) \in \bP(\cH)^2 \: \la x,y \ra \not=0\} \to \bP(\cH).\] 
\end{lem}

\begin{prf} In the following argument we represent elements of $\bP(\cH)$ by 
general non-zero vectors, not necessarily normalized to unit length. 
Recall that $\bP(\cH)$ is a symmetric space with the point reflections 
given by 
\[ r_{[x]}([y]) := \Big[-y + 2 \frac{\la y, x \ra}{\|x\|^2} x\Big].\] 
For the midpoint $[z] = [x] \sharp [y]$ of two non-orthogonal 
rays $[x]$ and $[y]$ we then have 
\[ [y] = r_{[z]}([x]) = \Big[- x + 2 \frac{\la x, z \ra}{\|z\|^2} z\Big].\]  
For $y = x + v$ and $z = x + w$ with $v, w\in x^\bot$ 
(here we normalize by $\la y,x \ra = \la z,x \ra = 1$), we then 
have 
\begin{align*}
[x + v] 
&= [y] = \Big[-x + \frac{2}{\|x + w\|^2} (x+w)\Big]
= \Big[-x + \frac{2}{1 + \|w\|^2} (x + w)\Big]\\
&= \Big[\frac{1-\|w\|^2}{1 + \|w\|^2}  x + \frac{2}{1 + \|w\|^2} w\Big] 
= \Big[x + \frac{2}{1 - \|w\|^2} w\Big].
\end{align*}
Note that this calculation actually shows that $\|w\| \not=1$ because 
$[y] \not=[w]$. We further obtain 
\[ v = \frac{2}{1 - \|w\|^2} w \quad \mbox{ and } \quad 
\|v\| = \frac{2\|w\|}{1 - \|w\|^2},\] 
which in turn yields 
\[ \|w\| = \frac{\sqrt{1 + \|v\|^2} - 1}{\|v\|}\] 
and hence 
\[ w = \frac{\sqrt{1 + \|v\|^2} - 1}{\|v\|^2}v.\] 
This argument also implies that $w$ is uniquely determined by 
$x$ and $y$ and that the midpoint operation is continuous in each argument separately. 

We now show that it is continuous. Pick  $[x_0] \in \bP(\cH)$ and let $U$ be 
an open neighborhood of $[x_0]$ for which there exists a continuous map 
$\sigma \: U \to \U(\cH)$ with $\sigma_u[x_0]= u$ for $u \in U$. 
For $\la x_0, y_0 \ra \not=0$ and 
$x_i \to x_0$, $y_i \to y_0$ we may w.l.o.g.\ assume that 
$\la x_i, y_i \ra \not=0$ for every $i$. We then have for 
$[x_i] \in U$ the relation 
\[ [x_i] \sharp [y_i]  
= \sigma_{[x_i]}([x_0] \sharp \sigma_{[x_i]}^{-1}[y_i])
\to [x_0] \sharp [y_0] \] 
because the action of $\U(\cH)$ on $\bP(\cH)$ is continuous. 
\end{prf}

\begin{lem} \mlabel{lem:pc} 
The set $\bP(\cH)_c$ has the following properties: 
  \begin{description}
  \item[\rm(i)] $\bP(\cH)_c$ is closed. 
  \item[\rm(ii)] If $[x], [y] \in \bP(\cH)_c$ are not orthogonal, then 
$[x,y] := \bP(\C x + \C y) \subeq \bP(\cH)_c$. 
  \item[\rm(iii)] If $[x_0], \ldots, [x_n] \in \bP(\cH)_c$ are such that 
$\la x_j, x_{j+1}\ra \not=0$ for $j =0,\ldots, n-1$, then 
$[x_0,\ldots, x_n] := \bP(\Spann\{x_0,\ldots, x_n\})\subeq \bP(\cH)_c$. 
  \end{description}
\end{lem}

\begin{prf} (i) follows immediately from Lemma~\ref{lem:closed}. 

(ii) We obtain from Lemma~\ref{lem:mipo} and the equivariance of the midpoint 
operation under $\PU(\cH)$ that, for two non-orthogonal elements $[x],[y] \in \bP(\cH)_c$, 
we also have $[x] \sharp [y] \in \bP(\cH)_c$. 

Let 
\[ \Exp \: T(\bP(\cH)) \to \bP(\cH) \] 
denote the exponential map of the symmetric space $\bP(\cH)$. 
For $[y] = \Exp(v)$ and  $v \in T_{[x]}(\bP(\cH))$ with $\|v\| < \pi/2$, 
the whole geodesic arc 
$[y_t] = \Exp(tv)$, $0 \leq t \leq 1$, consists of elements not orthogonal 
to $[x]$. By successive dyadic division, we can generate a dense subset 
of this arc by the midpoint operation from $[x]$ and $[y]$. 
Therefore $[x], [y] \in \bP(\cH)_c$ implies that 
$[y_t] \in \bP(\cH)_c$ for $0 \leq t \leq 1$. 

Since $\Exp$ is $\U(\cH)$-equivariant, we 
conclude that, for the action of $G$ on the tangent bundle 
$T(\bP(\cH))$, the set 
$T(\bP(\cH))_c$ of $G$-continuous elements has the 
property that, if $v \in T(\bP(\cH))_c$ with $\|v\| < \pi/2$, 
then $[0,1]v \subeq T(\bP(\cH))_c$. 
Since $G$ acts on $T(\bP(\cH))$ by bundle automorphisms, 
it also follows that, for each $[x] \in \bP(\cH)_c$, 
the set $T_{[x]}(\bP(\cH))_c$ is a closed complex linear subspace. 
This implies that, for two non-orthogonal elements 
$[x], [y] \in \bP(\cH)_c$, the whole projective 
plane $[x,y] \subeq \bP(\cH)$ consists of $G$-continuous vectors. 

(iii) We argue by induction and observe that the case $n = 1$ follows from (ii). 
Assume that $n > 1$. Then the induction hypothesis implies that 
$[x_0,\ldots, x_{n-1}] \subeq \bP(\cH)_c$. 
Since $x_n$ is not orthogonal to this space, 
the set of all elements $[y] \in [x_0,\ldots, x_{n-1}]$ with
$\la y,x_n \ra \not=0$ is open dense. Hence the closedness of 
$\bP(\cH)_c$ implies that $[x_0,\ldots, x_n] \subeq \bP(\cH)_c$. 
\end{prf}

\begin{defn} We call a subset 
$E \subeq \bP(\cH)$ {\it indecomposable} 
if 
\[ E = E_1 \dot\cup E_2, \quad E_1 \bot E_2, \quad 
E_1 \not=\eset \quad \Rarrow \quad E_2 = \eset,\] 
i.e., $E$ cannot be decomposed into two proper mutually orthogonal 
subsets. 
\end{defn}

\begin{rem} On every subset $E \subeq \bP(\cH)$ we obtain an equivalence 
relation by $[x] \sim [y]$ if there exists a sequence 
$[x_0], \ldots, [x_n] \in E$ with $[x_0] = [x]$, $[x_n] = y$ and 
$\la x_j, x_{j+1}\ra \not=0$ for $j =0,\ldots, n-1$. 
The corresponding equivalence classes are the maximal indecomposable subsets 
of $E$. We call them the {\it indecomposable components of $E$}. 
\end{rem}

\begin{thm} \mlabel{thm:pc} {\rm(Structure theorem for $\bP(\cH)_c$)} 
The indecomposable components of $\bP(\cH)_c$ are of the form 
$\bP(\cH_j)$, $j \in J$, where 
$(\cH_j)_{j \in J}$ is a family of mutually orthogonal closed subspaces of 
$\cH$. In particular, 
\[  \bP(\cH)_c = \bigcup_{j \in J} \bP(\cH_j),\] 
where each $\bP(\cH_j)$ is closed and open in $\bP(\cH)_c$. 
\end{thm}

\begin{prf} Let $C \subeq \bP(\cH)_c$ be an indecomposable component. 
Then $\bP(\cH)_c \subeq C \dot\cup C^\bot$ and the closedness of $C^\bot$ already 
implies that $C$ is relatively open in $\bP(\cH)_c$. Since this is also true for the 
other indecomposable components, $C$ is also relatively closed, 
hence closed in $\bP(\cH)$ because $\bP(\cH)_c$ is closed. 

For $[x], [y] \in C$ we find $[x_0],\ldots, [x_n] \in C$ with 
$[x] =[x_0]$, $[y] = [x_n]$ and 
$\la x_j, x_{j+1} \ra \not=0$ for $j =0,\ldots, n-1$. 
Therefore Lemma~\ref{lem:pc}(iii) implies that 
\[ [x,y] \subeq [x_0,\ldots, x_n] \subeq \bP(\cH)_c, \] 
and since $[x,y]$ is indecomposable, $[x,y] \in C$. 
This proves that $C = \bP(\cK)$ for a linear subspace 
$\cK \subeq \cH$. The closedness of $C$ now implies that $\cK$ is closed in $\cH$. 

If $\cP(\cK_1)$ and $\cP(\cK_2)$ are different indecomposable components of 
$\cP(\cH)_c$, then obviously $\cK_1 \bot \cK_2$. This completes the proof.
\end{prf}

\begin{cor} \mlabel{cor:cont} If $\bP(\cH)_c$ is total and indecomposable, 
then $\bP(\cH)_c = \bP(\cH)$. 
\end{cor}

\begin{cor} \mlabel{cor:cont2} If 
$\cP(\cK) \subeq \cP(\cH)_c$ is an indecomposable component, then 
\[ G_\cK := \{ g \in G \: \pi(g)\cP(\cK) \subeq \cP(\cK)\} \] 
is an open subgroup of $G$ and the induced homomorphism $G_\cK \to \PU(\cK)$ 
is continuous. 
\end{cor}

\begin{prf} Let $[v] \in \cP(\cK)$. Since $\cP(\cK)$ is an open subset 
of the $G$-invariant subset $\cP(\cH)_c$, the subset 
$U := \{ g \in G \: \pi(g)[v] \in \cP(\cK)\}$ is open. 
Since $G$ also permutes the indecomposable components of $\bP(\cH)_c$, 
any $g \in U$ preserves $\bP(\cK)$, so that $U = G_\cK$. 
The continuity of the projective representation of $G_\cK$ on $\cK$ 
follows from Corollary~\ref{cor:3.2}. 
\end{prf}

\begin{ex} Let $\chi \: \R \to \T$ be a discontinuous character 
and $\cH = \cH_1 \oplus \cH_2$ be a direct sum of two Hilbert spaces. 
Then $\pi(t)(v_1 + v_2) := v_1 + \chi(t) v_2$ for $v_j \in \cH_j$ 
defines a unitary representation on $\cH$ with 
$\cH_c = \cH_1$, but for the corresponding projective representation we have 
$\bP(\cH)_c = \bP(\cH_1) \dot\cup \bP(\cH_2)$. 
\end{ex}

\begin{thm} {\rm(Continuity criterion for projective representations)} 
\mlabel{thm:3.5} Let $G$ be a connected topological group and 
$\pi \: G \to \PU(\cH)$ be a projective unitary 
representation with cyclic ray $[v] \in \bP(\cH)_c$. Then $\pi$ is continuous. 
\end{thm}

\begin{prf} Since open subgroups of topological groups are also closed, 
the connectedness of $G$ implies that it preserves all indecomposable components 
of $\bP(\cH)_c$ (Corollary~\ref{cor:cont2}). Hence $\pi(G)[v]$ lies in a single indecomposable component 
$\bP(\cK)$ for a closed subspace $\cK \subeq \cH$. 
As $[v]$ is cyclic, we find that $\cK = \cH$, hence $\cP(\cH)_c = \cP(\cH)$, 
so that $\pi$ is continuous (Corollary~\ref{cor:3.2}). 
\end{prf}

\begin{cor} \mlabel{cor:3.7} Let $G$ be a connected topological group and 
$\pi \: G \to \PU(\cH)$ be an irreducible projective unitary 
representation with $\bP(\cH)_c\not=\eset$. Then $\pi$ is continuous. 
\end{cor}

\begin{thm} \mlabel{thm:4.17} Let $G$ and $T$  be topological groups 
and $\alpha \: T \to \Aut(G)$ be a homomorphims defining a continuous 
action of $T$ on $G$. Suppose that $G$ is connected 
and that $(\pi, \cH)$ is an irreducible continuous unitary representation of 
$G$ for which $\pi \circ \alpha_t$ is equivalent to $\pi$ for every $t \in T$. 
We then obtain a well-defined 
projective unitary representation 
$\oline U \: T \to \PU(\cH)$, defined by 
\[ U_t \pi(g) U_t^* = \pi(\alpha_t(g)) \quad \mbox{ for } \quad g \in G, t \in T,\] 
where $U_t \in \U(\cH)$ is a unitary lift of $\oline U_t \in \PU(\cH)$. 
Then the following assertions hold: 
\begin{description}
\item[\rm(i)] $\oline U$ is continuous if there exists a ray 
$[v] \in \PU(\cH)$ with a continuous orbit map. 
\item[\rm(ii)] If 
$\{ \ell_{\pi(g)v} \: g \in G \}$ separates the points of $\Herm_1(\cH)$, 
then $[v]$ has a continuous orbit map under $\oline U$ 
if and only if all function 
\[ T \to \R, \quad t \mapsto |\la U_tv, \pi(g)v \ra|, \quad g \in G,\] 
are continuous. 
 \item[\rm(iii)] If $\oline U$ is continuous, then 
\[ \tilde T  := \{ (t,U) \in T \times \U(\cH)_s \: 
\oline U_t  = \oline U\}, \quad p \: \tilde T \to T, 
(t,U) \mapsto t,  \] 
is a central $\T$-extension of $T$. 
By $\tilde\alpha(t,U) := \alpha(t)$ we then obtain a
 topological group $\tilde G := G \rtimes_{\tilde\alpha} \tilde T$ 
and $\tilde\pi(g,(t,U)) := \pi(g)U$ defines a continuous 
unitary representation of $\tilde G$ on $\cH$ extending~$\pi$. 
 \item[\rm(iv)] If, in addition, the central extension 
$p \: \tilde T \to T$ splits, then $\pi$ extends to a continuous unitary representation 
of $G^\sharp = G \rtimes_\alpha T$. This is always the case for $T = \R$. 
\end{description}
\end{thm}

\begin{prf} (i) Let $\bP(\cH)_c$ 
denote the set of continuous rays for the projective representation 
$\oline U$. Since $T$ acts continuously on $G$ and 
$\pi$ is continuous, the set $\bP(\cH)_c$ is $G$-invariant because 
\[ \oline U_t [\pi(g)v] = (U_t \pi(g) U_t^*) \oline U_t[v] 
= \pi(\alpha_t(g)) \oline U_t[v]. \] 

Suppose that $[v] \in \bP(\cH)_c$. Then 
$\pi(G)[v] \subeq \bP(\cH)_c$ is total and indecomposable (because $G$ is connected), 
so that Corollary~\ref{cor:cont} implies that $\bP(\cH)_c = \bP(\cH)$, i.e., 
that $\oline U$ is continuous. 

(ii) follows from Proposition~\ref{prop:2.9}. 
Note that $|\la U_tv, \pi(g)v \ra|$ depends only on 
$\oline U_t$ and not on the choice of the unitary lift $U_t$.

(iii) Since $\tilde T$ is the pullback $\oline U^*\U(\cH)$ of the 
central $\T$-extension $\U(\cH) \to \PU(\cH)$, it also is a 
central extension of topological groups, i.e., $p \: \tilde T \to T$ 
is a quotient map with continuous local sections. The rest of (iii) is 
clear.

(iv) If $\sigma \: T \to \tilde T, t \mapsto (t, U_t)$  is 
a continuous splitting of $\tilde T$, then 
$\pi^\sharp(g,t) := \pi(g) U_t$ is a continuous extension of  
$\pi$ to $G^\sharp$. 

For $T = \R$, the discussion in 
\cite[Ch.~9]{Mag92} implies that 
$\tilde T$ actually is a Lie group and since one-parameter groups 
of quotients of Lie groups lift, there exists a continuous splitting 
$\sigma \: T \to \tilde T$. 
\end{prf}

\section{Extending representations to semidirect products} 
\mlabel{sec:6} 

Let $G$ be a topological group. Then there exists a $W^*$-algebra 
$W^*(G)$, together with a homomorphism 
$\eta \: G \to \U(W^*(G))$ which is continuous with respect to the weak topology on 
$W^*(G)$, and which has the universal property that, for every continuous 
unitary representation $(\pi,\cH)$ of $G$, there exists a unique 
normal representation $\tilde\pi \: W^*(G) \to B(\cH)$ with 
$\tilde\pi \circ \eta = \pi$ (cf.\ \cite[Rem.~IV.1.2]{Ne00}, \cite{GN00}). 
In the appendix of \cite{GN00} it is shown that the predual 
$W^*(G)_*$ of this $W^*$-algebra can be identified 
with the subspace $B(G) \subeq C_b(G)$ of bounded continuous functions on $G$, 
spanned by the convex cone 
$\cP_c(G)$ of continuous positive definite function on $G$. 
The natural map 
\[ \eta_* \:  W^*(G)_* \to B(G), \quad \eta_*(\phi) := \phi \circ \eta \] 
is a linear bijection which is continous with respect to  
$\|\cdot\|_\infty$ on $B(G)$. 

\begin{rem} \mlabel{rem:5.1} In general the subspace $B(G)$ is not closed in $(C_b(G), \|\cdot\|_\infty)$, 
so that $\eta_*$ is not an open map with respect to $\|\cdot\|_\infty$ on $B(G)$. 

In fact, if $G$ is locally compact abelian, 
Bochner's Theorem implies that 
$B(G)$ is the range of the Fourier transform $M(\hat G) \to C(G)$ 
from the convolution algebra $M(\hat G)$ of finite Radon measures on $\hat G$ to $C(G)$.
For $G = \T$ and $\hat G = \Z$, we have $M(\hat G) \cong \ell^1(\Z)$, so that 
$B(\T)$ is the Wiener algebra of all continuous functions with absolutely convergent 
Fourier series. This algebra is dense in $C(\T)$ and a proper subspace, hence not closed. 
\end{rem} 

We endow $B(G)$ with the norm $\|\cdot\|$ that turns $\alpha_*$ into an isometry. 
According to \cite[Prop.~A.3]{GN00}, we then have 
$\|\phi\| = \|\phi\|_\infty = \phi(\1)$ for $\phi \in \cP_c(G)$,
but Remark~\ref{rem:5.1} implies that, in general, $\|\cdot\|$ is not equivalent 
to $\|\cdot\|_\infty$ because $B(G)$ may be incomplete w.r.t.\ $\|\cdot\|_\infty$. 

Let $T$ be a topological group and $\alpha \: T \to \Aut(G), t \mapsto \alpha_t$ 
be a homomorphism defining a continuous action of $T$ on $G$. We 
write $G^\sharp := G \rtimes_\alpha T$ for the corresponding semidirect product group. 
From the universal property of $W^*(G)$ we obtain a homomorphism 
$\tilde\alpha \: T \to \Aut(W^*(G))$ which is uniquely determined by 
$\tilde\alpha_t \circ \eta = \eta \circ \alpha_t$ for every $t\in T$. 
This defines a {\it $W^*$-dynamical system} in the sense of Borchers 
(where no continuity of $\tilde\alpha$ is required; cf.\ \cite{Bo83}).
Let $B(G)_c \subeq B(G)$ denote the set of all functions $\phi$ 
for which the map $T \to B(G), t \mapsto \alpha_t^*\phi 
= \phi \circ \alpha_t$ is continuous. 
From \cite[Prop.~II.5]{Bo83} it follows that $B(G)_c$ is a closed subspace
invariant under $T$ which is generated by the convex cone $B(G) \cap \cP_c(G)$. 

\begin{thm} \mlabel{thm:bo1} {\rm(Characterization Theorem)} 
Let $(\pi, \cH)$ be a continuous unitary representation of 
$G$ and $F_\pi \subeq B(G)$ the corresponding folium, i.e., the set of all 
functions of the form $\phi_S(g) := \tr(\pi(g)S)$, where 
$S \in \Herm_1(\cH)$ is non-negative with $\tr(S) = 1$. 
Then the following are equivalent: 
\begin{description}
\item[\rm(i)] $\pi$ is quasi-equivalent to a representation 
$(\pi',\cH')$ that extends to a continuous unitary representation 
$(\pi^\sharp,\cH')$ of $G^\sharp = G \rtimes_\alpha T$. 
\item[\rm(ii)] $F_\pi$ is $T$-invariant and contained in $B(G)_c$ 
(which implies that $T$ acts continuously on $F_\pi$). 
\end{description}
\end{thm}

\begin{prf} This follows immediately from \cite[Thm.~III.2]{Bo83}, applied to 
the $W^*$-dynamical defined by $(\tilde\alpha_t)_{t \in T}$ 
and the one-to-one correspondence between normal representations of 
$W^*(G)$ and continuous unitary representations of $G$. 
\end{prf}

\begin{cor} \mlabel{cor:5.3} {\rm(Extendability Criterion)} 
For a continuous positive definite function 
$\phi \in \cP_c(G)$, the following are equivalent:  
\begin{description}
\item[\rm(i)] $\phi$ extends to a continuous positive definite function 
of $G^\sharp = G \rtimes_\alpha T$. 
\item[\rm(ii)] $\phi \in B(G)_c$, i.e., the $T$-orbit map of $\phi$ is continuous with respect 
to the norm $\|\cdot\|$ on $B(G) \cong W^*(G)_*$. 
\end{description}
\end{cor}

\begin{prf} (i) $\Rarrow$ (ii): If $\phi^\sharp \: G^\sharp \to \C$ is a continuous 
positive definite extension of $\phi$, then the corresponding GNS representation 
of $G^\sharp$ is continuous, so that (ii) follows from Theorem~\ref{thm:bo1}. 

(ii) $\Rarrow$ (i): According to \cite{Bo93}, the subspace 
\[ \cM := \{ A \in W^*(G) \: A B(G)_c \cup B(G)_c A \subeq B(G)_c\}\] 
is a $W^*$-subalgebra, $\cN := B(G)_c^\bot \cap \cM$ is a $W^*$-ideal, 
and $B(G)_c$ is the predual of the $W^*$-algebra~$\cM/\cN$. 

Since the left and right multiplications with elements 
of $\eta(G) \subeq W^*(G)$ define isometries of $B(G) \cong W^*(G)_*$, 
and Theorem~\ref{thm:bo1} implies in particular that for each 
$\phi \in B(G)$ the maps 
\[ G \to B(G), \quad g\mapsto \eta(g)\phi, \quad 
g \mapsto \phi \eta(g) \] 
are continuous, Lemma~\ref{lem:closed} shows that the 
left and right multiplication actions of $G$ on $B(G)$ are continuous. 
From 
\[ \alpha_t^*(\eta(g)\phi) = \eta(\alpha_{t^{-1}} g) \alpha_t^*\phi \] 
we now conclude that $\eta(G) \subeq \cM$. Since 
$W^*(G)$ is generated by $\eta(G)$ as a $W^*$-algebra, 
it follows that $\cM = W^*(G)$. 
Any faithful normal representation of $\cM/\cN$ 
now yields a continuous unitary representation $(\pi, \cH)$ 
of $G$ for which 
\[ F_\pi = B(G)_c \cap \{ \phi \in \cP_c(G) \: \phi(\1) = 1\} \]  
is the corresponding folium. 
Since it is $T$-invariant, Theorem~\ref{thm:bo1} implies the 
existence of a quasi-equivalent representation $(\pi', \cH')$ of $G$ 
that extends to $G^\sharp$. As $F_\pi = F_{\pi'}$ by quasi-equivalence, 
(i) follows. 
\end{prf}

\begin{rem} \mlabel{rem:5.4} (a) 
Unfortunately, the characterization of the extendable positive definite functions 
on $G$ is given in terms of the norm $\|\cdot\|$ on $B(G)$ which is not very accessible. 
Since the inclusion $(B(G),\|\cdot\|) \to (C_b(G),\|\cdot\|_\infty)$ is continuous, 
Corollary~\ref{cor:5.3}  shows that the continuity of the map 
\[ T \to C_b(G), \quad t \mapsto \phi \circ \alpha_t \] 
is necessary for the existence of a continuous positive definite extension to $G^\sharp$. 

Since this condition is much easier to check for concrete cases,
it would be interesting to know whether it is also sufficient. 

(b) Assume that $T = \R$. 
Although, in general, the group $G^\sharp = G \rtimes_\alpha \R$ is not locally compact, 
the existence of an invariant measure on 
the quotient  $G^\sharp/G \cong \R$ implies that unitary induction makes sense 
as a passage from continuous unitary representations of $G$ to 
continuous unitary representations of $G^\sharp$. 

Starting with a unitary representation $(\pi, \cH)$ of $G$, we consider the space 
$\cH^\sharp := L^2(\R,\cH),$
endowed with the continuous unitary representation 
$\pi^\sharp := \Int_G^{G^\sharp}(\pi)$, given by 
\[ (\pi^\sharp(g,t)f)(s) = \pi(\alpha_{-s}(g))f(s-t) \quad \mbox{ for } \quad g \in G,s,t \in \R.\] 
Then $(\pi^\sharp(g,0)f)(s) = \pi(\alpha_{-s}(g))f(s)$ 
shows that $G$ acts by multiplication operators and 
$(U_tf)(s) := (\pi^\sharp(\1,t)f)(s) = f(s-t)$ shows 
that $\R$ acts by translations. It follows in particular that 
$\Spec(U) =\R$, so that $\pi^\sharp$ never is a positive energy representation. 

Let $v \in \cH$ be such that the positive definite function 
$\phi(g) := \la \pi(g)v,v\ra$ has the property that 
$\R \to C_b(G), t \mapsto \phi \circ \alpha_t$ is continuous. 
For $h \in C_c(\R)$, we consider the element $f(t) := h(t)v$ of $\cH^\sharp$. 
In the representation $(\pi^\sharp, \cH^\sharp)$ we then obtain the matrix coefficient 
\begin{align*}
 \la \pi^\sharp(g,t)f, f \ra 
&= \int_\R \la \pi(\alpha_{-s}(g)) v, v \ra h(s-t)\oline{h(s)}\, ds 
= \int_\R \phi(\alpha_{-s}(g)) h(s-t)\oline{h(s)}\, ds.
\end{align*}
This is a continuous positive definite function whose restriction to $G$ is given by 
\[  \la \pi^\sharp(g,0)f, f \ra 
= \int_\R \phi(\alpha_{-s}(g)) |h(s)|^2\, ds.\] 
If the functions $|h_n|_2^2$, $n \in \N$, form an approximate identity 
on $\R$, we thus obtain a sequence $(\phi_n^\sharp)_{n\in \N}$ 
of continuous positive definite functions on $G^\sharp$ such 
that $\phi_n := \phi_n^\sharp \res_{G}$ converges in $C_b(G)$ to $\phi$ 
(cf.\ \cite[Sect.~III]{Bo69}). 

If we apply the same construction to the discrete group $\R_d$ instead, we have 
$\cH^\sharp = \ell^2(\R,\cH)$ and for $f = \delta_0 v$, 
we obtain the positive definite function 
\[ \phi^\sharp(g,t) = \delta_{0,t} \phi(g)\]
which is not continuous on $G^\sharp$ if $\phi$ is not constant. 
\end{rem}

\subsection*{Applications to $C^*$-dynamical systems} 

We conclude this section with a brief discussion of the link 
to Borchers' criterion for the existence of covariant 
representations in the context of $C^*$-algebras. 

Let $\cA$ be a $C^*$-algebra, $T$ be a topological group and 
$\alpha \: T \to \Aut(\cA)$ be a group homomorphism 
(not necessarily strongly continuous). 
A covariant representation of $(\cA,\alpha)$ is a triple  
$(\pi, U, \cH)$, where 
$\pi \: \cA \to B(\cH)$ is a non-degenerate representation of $\cA$  
and $U  \: T \to \U(\cH)$ is a continuous representation satisfying 
\[ U_t \pi(A) U_t^*  = \pi(\alpha_t(A)) \quad \mbox{ for } \quad 
A \in \cA, t \in T.\] 

If $(\pi, U,\cH)$ is a continuous unitary representation 
and $v \in \cH$ a unit vector, then 
$\phi(A) := \la \pi(A)v,v\ra$ is a normalized state 
$\phi \in \cA^*$ for which the $\alpha$-orbit map 
$T \to \cA^*, t \mapsto \alpha_t^*\phi := \phi \circ \alpha_{t}$ is norm continuous. 
This follows immediately from 
\begin{equation} \label{eq:covrelx}
(\alpha_t^* \phi)(A) 
= \la \pi(\alpha_t A)v,v\ra
= \la U_{t} A U_{t}^*v,v\ra 
= \la A U_t^*v, U_t^*v\ra= \tr(A P_{U_t^*v}).
\end{equation}
There is an interesting converse result, namely that 
every state $\phi \in \cA^*$ with a norm-continuous 
$\alpha$-orbit map comes from a vector in a covariant 
representation (cf.\ \cite{Bo83}). Since Borchers' proof of this result 
is quite involved, one would like to a have a more direct argument 
based on the criteria from above.

\begin{rem} \mlabel{rem:4.8} (cf.\ \cite{BN12}) 
Let $\pi \: \cA \to B(\cH)$ be a representation of the 
$C^*$-algebra and 
\[ \pi_* \: B_1(\cH) \to \cA^*,\quad \pi_*(X)(A) := \tr(\pi(A)X).\] 
Then $\pi_*$ is a continuous linear map whose adjoint 
$\pi_*^* \: \cA^{**}  \to B_1(\cH)^* \cong B(\cH)$ 
is a normal representation of the $W^*$-algebra $\cA^{**}$ 
(cf.\ \cite{Sa71}). 
In particular, its range is closed, and therefore the range of 
$\pi_*$ is closed as well. Therefore $\pi_*$ induces a topological 
embedding $B_1(\cH)/\ker \pi_* = B_1(\cH)/\pi(\cA)^\bot \into \cA^*$ 
of Banach spaces. 

The map $\pi_*$ is injective if and only if $\pi(\cA)$ is weakly dense in 
$B(\cH)$, which is equivalent to $\pi(\cA)'' = \pi_*^*(\cA^{**}) = B(\cH)$. This 
means that $\pi$ is irreducible and then $\pi_* \: B_1(\cH) \to \cA^*$ is a topological 
embedding. 
\end{rem}

\begin{thm} \mlabel{thm:4.19} Assume that $\alpha \: T \to \Aut(\cA)$ 
defines a strongly continuous automorphic $T$-action on~$\cA$. Let 
$(\pi, \cH)$ be an irreducible representation 
of $\cA$ for which 
$\pi \circ \alpha_t$ is equivalent to $\pi$ for every $t \in T$ 
and let $v \in \cH$ be a unit vector. Then the following are equivalent: 
  \begin{description}
    \item[\rm(i)] The state $\phi(A) := \la \pi(A)v,v\ra$ of $\cA$ has a 
norm-continuous orbit map $T \to \cA^*, t \mapsto \phi \circ \alpha_t$. 
    \item[\rm(ii)] There exists a central $\T$-extension 
$p \: \tilde T \to T$ and a continuous unitary representation 
$\tilde U \: \tilde T \to \U(\cH)$ 
such that $(\pi, \tilde U,\cH)$ is a covariant 
representation with respect to $\tilde \alpha := \alpha \circ p 
\: \tilde T \to \Aut(\cA)$. 
  \end{description}
\end{thm}

\begin{prf} In view of \eqref{eq:covrelx},  (ii) implies (i). 
So we assume (i). Then there exists a family 
$(U_t)_{t \in T}$ of unitary operators with  
\begin{equation}
  \label{eq:covar2}
 U_t \pi(A) U_t^* = \pi(\alpha_t(A)) \quad \mbox{ for } \quad t \in T, A \in \cA.
\end{equation}

From \eqref{eq:covrelx} we derive the relation 
$\alpha_t^* \phi = \pi_*(P_{U_{t^{-1}}v}).$
Since $\pi_*$ is a topological embedding by Remark~\ref{rem:4.8}, the map 
$T \to \Herm_1(\cH), t \mapsto P_{U_{t}v}$ is also continuous. 
Therefore Lemma~\ref{lem:2.1} shows that 
$[v] \in \bP(\cH)_c$ is a continuous ray for the corresponding 
homomorphism $\oline U \: T \to \PU(\cH)$. 
As in the proof of Theorem~\ref{thm:4.17}, we see that 
$\pi(\U(\cA))[v]\subeq \bP(\cH)_c$ consists of continuous rays for~$\oline U$. 
As $\pi(\U(\cA))$ acts transitive on $\U(\cH)$ 
(\cite[Thm.~2.8.3(iii)]{Dix64}), we obtain 
$\bP(\cH) = \bP(\cH)_c$, so that Corollary~\ref{cor:3.2} implies that 
the $T$-action on $\bP(\cH)$ is continuous 
and now one proceeds as in the proof of Theorem~\ref{thm:4.17}. 
\end{prf}

\begin{rem} 
If $\pi$ is not irreducible, then the situation is more complicated. 
Then $\pi(\cA)' \not=\C\1$, so that the ambiguity in the choice of $U$ is rather 
large and the requirement that $\pi \circ \alpha_t \sim \pi$ for every 
$t$ does not lead to a canonical projective unitary representation.
\end{rem}

\begin{ex}   Let $\cA$ be a $C^*$-algebra and 
$\alpha \: T \to \Aut(\cA)$ a homomorphism defining a 
continuous action of $T$ on $\cA$. 
Then $\alpha$ defines in particular a continuous $T$-action on the 
unitary group $G := \U(\cA)$, endowed with the norm topology. 
Let $\phi \: G \to \C$ be a positive definite function defined by a 
linear positive functional $\psi \in \cA^*$. 
We claim that, in this case, the requirement that 
the curve $\phi_t$ in $C_b(G)$ is  continuous 
implies that the curve $\psi_t := \psi \circ \alpha_t$ 
in $\cA^*$ is norm continuous 
(cf.\ Remark~{rem:5.4}(a)). 

This claim follows from the fact that 
the absolute convex hull of $\U(\cA)$ is a $0$-neighborhood in $\cA$. 
In fact, if $A = A^*$ with $\|A\| < 1$, then 
$U := A + i \sqrt{\1 - A^2}$ 
is unitary and $A = \shalf(U + U^*)$. 
Therefore the absolute convex hull of $\U(\cA)$ 
contains all hermitian operators $A$ with $\|A\| < 1$, 
hence all operators $C$ with $\|C\| < \shalf$ 
because $C = \shalf (C + C^*) -i \shalf (iC - i C^*)$. 
\end{ex}

\end{document}